\documentclass[10pt]{article}
\usepackage{amssymb,latexsym, amsmath}
\usepackage{amscd}
\newtheorem{remark}{\it Remark}
\newtheorem{theorem}{Theorem}
\newtheorem{proposition}{Proposition}
\newtheorem{corollary}{Corollary}
\newtheorem{definition}{Definition}
\newtheorem{conjecture}{Conjecture}
\newtheorem{lemma}{Lemma}

\newcommand{\St}{\mathrm{St\,}}
\newcommand{\Cal}{\mathcal}
\newcommand{\goth}{\mathfrak}

\DeclareMathOperator{\ad}{\mathrm{ad}}
\DeclareMathOperator{\Ad}{\mathrm{Ad}}

 \DeclareMathOperator{\pr}{\mathrm{pr}}
\DeclareMathOperator{\trdeg}{\mathrm{tr.deg}}

\DeclareMathOperator{\ind}{\mathrm{ind}}
\DeclareMathOperator{\Ann}{\mathrm{Ann}}
\DeclareMathOperator{\Frac}{\mathrm{Frac}}
\DeclareMathOperator{\Ker}{\mathrm{Ker}}

 \DeclareMathOperator{\codim}{\mathrm{codim}}
 \DeclareMathOperator{\Tr}{\mathrm{Tr}}

\begin{document}

\title{Complete commutative subalgebras
in polynomial Poisson algebras:  a proof of the Mischenko--Fomenko conjecture\footnote{This text is a revised version of my paper published in Russian: A.V.Bolsinov, Complete commutative families of polynomials in Poisson-Lie algebras: A proof of the Mischenko-Fomenko conjecture//In book: Tensor and Vector Analysis, Vol.~26, Moscow State University, 2005, pp.~87-109.  The present English version has been available on my home page {\tt http://www-staff.lboro.ac.uk/\,$\widetilde{}$ maab2/} since 2008.  I apologise that  I have not changed anything in the text since then,  and some references are unfortunately out-of-date.}}

\author
{Alexey V. Bolsinov}

\maketitle

\abstract{
The Mishchenko-Fomenko conjecture  says that for each real or complex
finite-dimensional Lie algebra $\goth g$  there exists a complete set
of commuting polynomials on its dual space $\goth g^*$. In terms of the
theory of integrable Hamiltonian systems this means that the dual space
$\goth g^*$ endowed with the standard Lie-Poisson bracket admits
polynomial integrable Hamiltonian systems. Recently this conjecture has
been proved by S.T.~Sadetov.  Following his idea, we give an explicit geometric
construction for commuting polynomials on $\goth g^*$ and consider some
examples. }

\section{
Introduction and preliminaries
}

Consider a symplectic manifold $(M^{2n},\omega)$ and a Hamiltonian system $\dot x=X_H(x)$ on it, where $H:M^{2n}\to \Bbb R$ is a smooth function called {\it Hamiltonian}  and  $X_H(x)=\omega^{-1}(dH(x))$ is the corresponding Hamiltonian vector field.

This system is called {\it completely integrable} if  it admits  $n$ functionally independent  integrals $f_1,\dots, f_n:M^{2n} \to \Bbb R^n$ which commute with respect to the Poisson bracket associated with the symplectic structure $\omega$,  i.e., $\{f_i,f_j\}=0$, $i,j=1,\dots,n$.

Equivalently one can say that this system admit a complete commutative subalgebra $\Cal A$ of integrals in the Poisson algebra $C^\infty (M^{2n})$ of smooth functions on $M$.
Completeness means that at a generic point $x\in M^{2n}$, the subspace in $T^*M$ generated by
the differentials $df(x)$, $f\in\Cal A$  is maximal isotropic.

The same definition makes sense if, instead of a symplectic manifold, we consider a Poisson manifold
$(M, \{,\})$ where the Poisson bracket $\{,\}$ is not necessarily non-degenerate.

One of the most intriguing questions in the theory of integrable systems can be formulated as follows:
does a given symplectic (Poisson) manifold $M$ admit an integrable system with nice properties?

Notice that the necessity of "nice properties" is motivated by the fact  that any symplectic (Poisson) manifold admits a smooth integrable system which can be constructed by using  some kind of "partition of unity" idea \cite{Fomsymp}.  The behavior of such a system, however,  has no relation to the geometry of the underlying manifold and therefore  is not of interest at all.

The additional assumptions that make the above question non-trivial and interesting can be rather
various.  Briefly, we mention three types of integrable systems for which the existence problem is extremely interesting and important:

1) toric (or almost toric) integrable systems
\cite{Delz, Aud, Sym, Vu};

2) integrable systems with non-degenerate singularities \cite{El}, \cite{Zung1}, \cite{Zung2}, \cite{BO};

3) integrable geodesic flows on compact manifolds \cite{Koz}, \cite{But}, \cite{PS}, \cite{BJbelg}.

In the algebraic case, the existence problem seems to be interesting even without any additional assumptions:  given an algebraic symplectic (Poisson) manifold $X$, does it admit
a polynomial (rational) integrable system? In the present  paper,   we discuss this problem in the case  when $X$ is a dual space of a finite-dimensional Lie algebra endowed with  the standard linear Lie-Poisson bracket.

We start with recalling basic definitions. Consider  a finite-dimensional Lie algebra $\goth g$ over $\Bbb R$
and its dual space $\goth g^*$ endowed with the standard
Poisson-Lie structure which is defined as follows.
Let $f,g:\goth g^*\to\Bbb
R$ be arbitrary smooth functions. Their differentials at a point
$x\in \goth g^*$ can be treated as elements of the Lie algebra $\goth g$.
Then the  Lie-Poisson bracket of $f$ and $g$
is defined by:
\begin{equation}
\{ f,g\}(x)= \langle x,[df(x),dg(x)]\rangle.
\label{bracket}
\end{equation}

%%Symmetric algebra  Notation!!!!!!  Change S(\goth g)   --->  S(\goth g)

If instead  of smooth functions  we
restrict ourselves with polynomials on $\goth g^*$, then the same operation can
be introduced in the following equivalent way. The Poisson-Lie bracket on
the space of polynomials is defined to be a bilinear skew-symmetric
operation satisfying two properties:

1)  $\{ fg, h\}= f\{g,h\} + g\{f,h\}$ (Leibniz rule);

2) if  $f,g\in \goth g$ are linear polynomials on $\goth g^*$
then the Poisson-Lie
bracket coincides with the usual commutator in $\goth g$,  i.e.,
$$ \{f,g\}=[f,g].  $$

The space of polynomials $\Bbb R[\goth g]$ with such an operation is called the
Poisson algebra (associated with $\goth g$) and is denoted by $S(\goth g)  $.

The Poisson-Lie bracket is naturally extended to the space of rational
functions
$\Bbb R(\goth g)=\Frac (S(\goth g)  )$, and (which is very important for our
considerations) all the definitions make sense over arbitrary field $\Bbb K$ of  zero
characteristic.

To each finite-dimensional Lie algebra (over a field $\Bbb K$) one can
assign two integer numbers:  its dimension $\dim \goth g$ and index $\ind \goth g$.
The latter is the corank of the skew-symmetric form $\Phi_x : \goth g\times \goth g
\to \Bbb K$ for a generic element $x\in \goth g^*$ where
$$
\Phi_x(\xi,\eta)=\langle x, [\xi,\eta] \rangle.
$$

\begin{definition}\label{compl}
A commutative set of algebraically independent polynomials
$$
f_1,\dots,f_k
\in S(\goth g)
$$
is called complete, if $k=\frac{1}{2}(\dim \goth g + \ind \goth g)$.

A commutative subalgebra $\Cal A\subset S(\goth g)$ is called complete if
$\trdeg \Cal A = \frac{1}{2}(\dim \goth g + \ind \goth g)$.
\end{definition}

The completeness condition means that, at a generic point
$x\in \goth g^*$, the subspace
in $\goth g$ generated by the differentials
$df_1(x),\dots,df_k(x)$ is maximal isotropic with respect to the
Lie-Poisson bracket at $x$, i.e., in the sense of the skew-symmetric form
$\Phi_x$. In particular, the maximal possible number
of commuting independent polynomials in $S(\goth g)  $ cannot exceed
$\frac{1}{2}(\dim \goth g + \ind \goth g)$.

\begin{conjecture} {\rm (Mishchenko-Fomenko \cite{MFtr})}
Let  $\goth g$ be a real or complex finite-dimensional Lie algebra.
Then on $\goth g^*$ there exists a complete commutative set of polynomials.
\end{conjecture}

In more algebraic terms this means that each Poisson algebra $S(\goth g)  $
admits a complete commutative subalgebra $\Cal A$.

This conjecture comes from the theory of integrable Hamiltonian systems
and can be reformulated as follows: on the dual space $\goth g^*$
of every finite-dimensional Lie algebra $\goth g$ there exist integrable
Hamiltonian systems with polynomial integrals.

In 1978 A.~Mishchenko and A.~Fomenko \cite{MFizv} proved this conjecture for semisimple
Lie algebras. Since then complete commutative sets have been constructed
for many other classes of Lie algebras (see \cite{FomTrof}, \cite{Bolsizv}, \cite{Trof}, \cite{Thimm}). Recently S.~Sadetov
\cite{Sad}
has proved this conjecture in the general case by using one nice algebraic
construction  that reduces the problem either to the
semisimple case, or to an algebra of smaller dimension.

\begin
{theorem} {\rm (Sadetov, 2003)}  The Mishchenko-Fomenko conjecture holds
for an arbitrary finite-dimensional Lie algebra over a field of
zero characteristic.
\label{theo1}
\end{theorem}

It is remarkable fact that working over an arbitrary field surprisingly simplifies the proof. The main construction is based on the induction argument.  On each step we reduce the dimension of the
Lie algebra in question, but we have to pay for this by extending the field.  However, this price is not very high since all the statements and definitions admit purely algebraic formulations so that the field do not play any essential role.

The purpose of this paper is to present Sadetov's construction in  more
explicit terms of Poisson geometry allowing one to work effectively with
specific Lie algebras.  The approach suggested by S.Sadetov  is, in fact, purely algebraic.
In our opinion, however, behind his construction one can see  important geometrical ideas which we
would like to emphasize in the present paper rather that to  give another rigorous proof.
We also study several natural examples of Lie algebras and describe explicitly the related complete commutative subalgebras some of which are quite remarkable.

The proof which we are going to present is actually based on a
modification  of two well-known constructions: the "argument shift" method
suggested by A.~Mishchenko and A.~Fomenko and the so-called "chain of
subalgebras" method which was used by many authors for different purposes
(see, in particular, Gelfand-Zetlin \cite{GZ}, Vergne \cite{Vergne}, Thimm \cite{Thimm}, Trofimov \cite{Trof}). We start with recalling these constructions.

%\newpage

\section{
"Chain of subalgebras" method
}

In this section, by $\goth g$ we mean a real or complex Lie algebra. However
almost all constructions make sense for any field
$\Bbb K$ of zero characteristic.

Let $\goth h\subset \goth g$ be a subalgebra. Suppose that we can construct a
complete commutative subalgebra $\Cal A$  in $S(\goth h) $. Since $S(\goth h)  \subset
S(\goth g)  $, we can try to extend $\Cal A$ up to a complete commutative subalgebra in $S(\goth g)  $. To this end we need to find additional polynomials $f_1,\dots, f_s$
which commute with $S(\goth h)  $ and between themselves. As good candidates we
can use, for examples, the invariants of the coadjoint representation of
$\goth g$ or, which is the same, the polynomials from the center $Z(S(\goth g)  )$ of
$S(\goth g)  $.  Sometimes these polynomials are sufficient to satisfy the
completeness condition.

Repeating this idea for a chain of subalgebras
$$
\{0\} = \goth g_0 \subset \goth g_1 \subset \goth g_2 \subset \dots \subset \goth g_{n-1} \subset
\goth g_{n} = \goth g
$$
we can always construct a "big" set of commuting polynomials:
$$
Z_0\cup Z_1 \cup \dots \cup Z_{n-1} \cup Z_{n},
$$
where $Z_i=Z(S(\goth g_i))$.

For many important cases this allows us  to construct a complete commutative subalgebra  in $S(\goth g)  $. For example, it is so for the
chains  {(see \cite{Thimm})}
$$
gl(1,\Bbb R) \subset gl(2,\Bbb R) \subset \dots \subset gl(n-1,\Bbb R)
\subset gl(n,\Bbb R),
$$
$$
so(1,\Bbb R) \subset so(2,\Bbb R) \subset \dots \subset so(n-1,\Bbb R)
\subset so(n,\Bbb R),
$$
and also for codimension one filtrations in
nilpotent (see \cite{Vergne}) and solvable algebraic Lie algebras (in the latter case
instead of polynomials one has to consider rational functions, but after
some modification using semi-invariants instead of invariants one still can
solve the problem without leaving the space of polynomials).

However in the general case an appropriate chain of subalgebras does
not always exist, and this method does not work directly.

Let us look at the problem with more attention.  To understand the situation better, let us
first consider the following "linear" version of our  problem.  Take a vector space $V$ endowed with a skew-symmetric bilinear form  $\phi$ (possibly, degenerate!).  Let $U_1 \subset V$ be a subspace, and $A_1\subset U_1$ be a maximal  isotropic subspace in $U_1$.   The problem is to extend $A_1$ up to a maximal isotropic subspace $A\subset V$.  One of  possible solutions is the following.  Consider the
skew-orthogonal "complement" of $U_1$ in $V$, i.e., subspace
$$
U_2= U_1^\phi=\{ v\in V~|~\phi(u,v)=0  \, \, \mbox{for all} \,\, u\in U_1\}.
$$
Let $A_2\subset U_2$ be a maximal isotropic subspace in $U_2$. Then $A=A_1+A_2$ is maximal isotropic in $V$. This is a simple fact from linear symplectic geometry.

We now consider a "non-linear" version of this statement.  Consider a Poisson manifold
$(X, \phi)$ and a (Poisson) subalgebra $\Cal F\subset C^\infty(X)$.  A commutative subalgebra $\Cal A\subset \Cal F$ is called complete in $\Cal F$ if  at a generic point $x\in M$ the following condition holds.
Consider the subspaces $d\Cal A(x)$ and $d\Cal F(x)$ in $T^*_xX$ generated by the differentials of functions $f$ from $\Cal A$ and $\Cal F$ respectively.  It is clear that  $d\Cal A(x)$ is an isotropic subspace in  $d\Cal F(x)$  with respect to the Poisson structure $\phi$.

\begin{definition}
A commutative subalgebra $\Cal A\subset \Cal F$ is called complete in $\Cal F$ if  $d\Cal A(x)$ is {\it maximal} isotropic in $d\Cal F(x)$  at a generic point $x\in X$.
\end{definition}

 Now consider two (Poisson) subalgebras $\Cal F_1, \Cal F_2\subset C^\infty (X)$ such that
$\{ \Cal F_1, \Cal F_2\}=0$.  Let $\Cal A_1\subset \Cal F_1$, $\Cal A_2\subset \Cal F_2$ be complete commutative subalgebras in $\Cal F_1$ and $\Cal F_2$ respectively.
The following proposition is just a reformulation of the "linear" statement.

\begin{proposition}\label{prop1}
Suppose $d\Cal F_2(x) = d\Cal F_1(x)^\phi = \{  \xi \in T^*_x X~|~  \phi(\xi, df(x))=0 \  \mbox{for any } \  f\in \Cal F_1\}$ at a generic point $x\in X$.  Then
$\Cal A_1 + \Cal A_2$ is a complete commutative subalgebra in $C^\infty(X)$.
\end{proposition}

Here by  {\it generic} we mean "from open everywhere dense subset" without specifying the nature of such a subset,    and $\Cal A_1 + \Cal A_2$ denotes the least Poisson subalgebra in $C^\infty(X)$ which contains  both $\Cal A_1$ and  $\Cal A_2$.

\begin{remark}
{\rm The condition $d\Cal F_2(x) = d\Cal F_1(x)^\phi$ can be replaced by the following assumption:
$d\Cal F_2(x) + d\Cal F_1(x)$ is coisotropic in $T^*_x(X)$, which is slightly weaker.}
\label{rem1}
\end{remark}

This simple idea can now be applied to our problem.
Having a complete commutative subalgebra $\Cal A\subset S(\goth h)  $, we need to extend it up to
a complete commutative subalgebra in $S(\goth g)  $. Following the above construction, we should consider the maximal subalgebra  in $S(\goth g)  $  all of whose elements commute with $S(\goth h)  $.  Since $S(\goth h)  $ is generated by $\goth h$,  this subalgebra is:
$$
\Ann(\goth h)=\{ f\in S(\goth g)  ~|~ \{f, \eta\}=0, \quad \forall \ \eta \in \goth h \}.
$$

It is easy to see that $\Ann(\goth h)$  consists exactly of invariant polynomials
with respect to the coadjoint action of $H$ on $\goth g^*$,
where $H\subset  G$ is the Lie subgroup corresponding to $\goth h$.
To apply Proposition \ref{prop1}  we have to assume that this representation admits sufficiently many polynomial invariants.  More precisely, this means that elements of
$\Ann(\goth h)$  distinguish generic orbits, i.e.,
\begin{equation}
 \trdeg \Ann(\goth h) =\codim \Cal O_H(x),
\label{cond1}
\end{equation}
where  $\Cal O_H(x)\subset \goth g$ is a generic $\Ad^*_H$-orbit.

Notice that this condition means exactly that
$$
d\Ann(\goth h)(x) =  \goth h^{\Phi_x}=\{
\xi\in \goth g~|~ \langle x,  [\xi,\eta] \rangle=0 \ \forall \eta\in \goth h \}
$$
at a generic point $x\in \goth g^*$ and we can reformulate Proposition \ref{prop1} as follows.

\begin{proposition} Let $\goth h$ and $\Ann(\goth h)$ both admit complete
commutative subalgebras of polynomials $\Cal A_1\subset S(\goth h)  $ and $\Cal A_2\subset \Ann(\goth h)$
respectively.  If  {\rm (\ref{cond1})} holds, then $\Cal A_1 + \Cal A_2$ is a complete commutative subalgebra in $S(\goth g)  $.
\label{prop2}
\end{proposition}

\begin{remark}
{\rm
If we work over an arbitrary field $\Bbb K$ of zero characteristic,  then condition (\ref{cond1}) is not so convenient and can be replaced by one of the two following assumptions which do not involve any Lie groups:

1) $d\Ann(\goth h)(x) + \goth h$
is a coisotropic subspace in $\goth g$  w.r.t. $\Phi_x$  for generic $x\in \goth g^*$.

2) $\trdeg \Ann(\goth h) = \codim  \ad^*_{\goth h} x$  for generic $x\in \goth g^*$  (in the classical case where  $\Bbb K = \Bbb C$ or $\Bbb R$,   this subspace
$\ad^*_{\goth h}x \subset \goth g^*$  is just the tangent space for the orbit $\Cal O_H(x)$ at $x$).
 }
\label{rem2}
\end{remark}

Thus,  to construct a complete commutative subalgebra in $S(\goth g)  $,
it suffices to find complete commutative subalgebras in $S(\goth h)  $ and
$\Ann(\goth h)$.
Usually the dimension of  $\goth h$ and the transcendence degree of $\Ann(\goth h)$ are both  smaller than $\dim \goth g$, and we may hope that the problem of constructing complete
commutative subalgebras in $S(\goth h)$ and $\Ann(\goth h)$ will be simpler than that in
$S(\goth g)  $. The difficulty, however, is that $\Ann(\goth h)$ may have a
rather complicated algebraic structure.

It appears, nevertheless, that each non-semisimple Lie algebra always admits an ideal
$\goth h\subset \goth g$  such that $\Ann(\goth h)$  has a very nice structure. Roughly speaking,
$\Ann(\goth h)$  can be treated as a symmetric algebra $S(L)$  of a certain finite-dimensional Lie algebra $L$  but perhaps over a new field
$\Bbb K$.  After this,
according to Proposition  \ref{prop2} our problem is reduced to the same problem for smaller algebras $\goth h$ and $L$, which allows us to use the induction argument.

%\newpage

\section{
Argument shift method
}

The argument shift method was suggested by A.T.~Fomenko and A.S.~Mishchenko
in \cite{MFizv} as a generalization of S.V.~Manakov's construction \cite{Man}.

Let $\goth g$ be a Lie algebra, $\goth g^*$ be its dual space. Consider the
ring of invariants of the coadjoint representation $\Ad^*:  G \to
GL(\goth g^*)$:
$$
I_{\Ad^*}(G)=\{ f:\goth g^* \to \Bbb R~|~ f(l)=f(\Ad^*_g l) \text{ \ \ for any
$g\in G$} \}
$$

Generally speaking, the $\Ad^*$-invariants are not necessarily
polynomials. But locally in a neighdorhood of a regular element $x\in \goth g^*$
we always can find $k=\ind \goth g$ functionally independent smooth invariants.

For a fixed regular element $a\in \goth g^*$, consider the family of
functions
$$ \Cal F_a = \left\{ f_\lambda(x)=f(x+\lambda a) \right\}_{f\in
I_{\Ad^*}(G), \lambda\in \Bbb R}.  $$

It turns out that this family is commutative with respect to the
Lie-Poisson structure.
As we already noticed, the commuting functions so obtained are not
necessarily polynomials. However, this trouble can be avoided by replacing
the functions $f(x+\lambda a)$ with the homogeneous polynomials
$f_k(x)$ obtained by Taylor expansion of $f(x)$ at the point $a\in \goth g^*$:
$$
f(a+\lambda x) = f(a) + \lambda f_1(x) + \lambda^2 f_2(x) + \dots
$$

As a result, we
shall obtain a commutative subset $\{ f_k \}_{f\in I_{\Ad^*}(G)}$ in
$S(\goth g)  $ which we shall still denote by $\Cal F_a$.

\begin{theorem} {\rm (Mishchenko, Fomenko \cite{MFizv})} If $\goth g$ is semisimple
and $a\in \goth g^*$ is a regular element, then the commutative set $\Cal
F_a$ is complete.
\label{MF}
\end{theorem}

It is well known that  the argument shift method is closely related to compatible Poisson brackets and bi-Hamiltonian systems. Indeed,  on $\goth g^*$ there are two natural compatible Poisson brackets. The first one is the standard Poisson-Lie bracket (\ref{bracket}), the second is given by
$$
\{ f,g\}_a(x)= a([df(x),dg(x)]),
$$
where $a\in \goth g^*$ is a fixed element.

The compatibility condition is straightforward and the bi-hamiltonian approach (see \cite{Magri}) leads us immediately to Hamiltonian systems whose first integrals are Casimir functions of  linear combinations $\{,\}+\lambda \{,\}_a$,  which coincide exactly with the functions from $\Cal F_a$.

The bi-hamiltonian approach can be applied for an arbitrary Lie algebra, not necessarily semi-simple, and in fact, the family  $\Cal F_a$  turns out to be complete  for many other classes of Lie algebras.
More precisely, the following criterion holds.

Consider the set of singular elements in $\goth g^*$:
$$
\mathrm{Sing}=\{ l\in \goth g^*~|~ \dim\St_{\ad^*}(l) > \ind \goth g\},
$$
where $\St_{\ad^*}(l)=\{ \xi\in\goth g~|~ \ad^*_\xi l=0 \}$ is the
stationary subalgebra of $l$ in the sense of the coadjoint representation.

If $\goth g$ is an algebra over $\Bbb R$, then $\mathrm{Sing}$ is taken in the
complexification $(\goth g^{\Bbb C})^*$.

\begin
{theorem} {\rm (\cite{Bolsizv})}  Let $a\in \goth g^*$ be a regular element. The commutative set
$\Cal F_a \subset S(\goth g)  $ is complete if and only if $\codim \mathrm{Sing} >1$.
\end{theorem}

It is important to remark that in the
semisimple case the argument shift method works for any field of zero
characteristic.  This follows from the fact that the completeness
condition is preserved under extension of the field.

We now consider an example of a semisimple Lie algebra over a
"non-standard" field to show how the argument shift
methods works in a more complicated situation.

Consider a linear representation $\rho$ of a complex Lie algebra $\goth g$ on a
linear space $V$.

Consider all rational mappings $\Psi : V \to \goth g$ satisfying the following
property:  $\Psi(v)\in \St(v)$ where $\St(v)
=\{ \xi\in \goth g~|~ \rho(\xi) v =0\}$ is the stationary subalgebra of $v$ with
respect to $\rho$.

In other terms, $\Psi$ can be treated as a rational section of the
stationary subalgebra fiber bundle over $V$ (the fact that these
subalgebras are of different dimensions is not important, over an
Zariski open set this fiber bundle is smooth and locally trivial).

It is easy to see that the space $L=L(\goth g,\rho, V)$ of such sections can be
endowed with a Lie algebra structure. Indeed, we can just put by
definition:  $$ [\Psi_1, \Psi_2] (v) = [\Psi_1(v),\Psi_2(v)] \in \St(v).
$$

Over the original field this Lie algebra
$L=L(\goth g,\rho, V)$ is infinite dimensional.  But, we can, obviously,
consider it over the field $\Bbb K=\Bbb C(v_1,\dots,v_k)$ of rational
functions on $V$. Then $L(\goth g,\rho, V)$ has a finite dimension and, moreover,
$\dim_{\Bbb K} L(\goth g,\rho, V)$ is equal to the dimension (over $\Bbb C$) of
a generic stationary subalgebra.

Assume that a generic stationary subalgebra $\St(v)$ is semisimple, then
so is $L(\goth g,\rho, V)$ over $\Bbb K$.

Let us construct a complete commutative set in $P(L(\goth g,\rho, V))$ by using
the argument shift method. First of all we notice that, as usual,
$L^*$ can be identified with $L$ (and, consequently,
$\ad$ with $\ad^*$) by using the form $\Tr: L\times L \to\Bbb K$:
$$
(\Tr \Psi_1\Psi_2) (v)=\Tr_\rho \bigl( \Psi_1(v)\Psi_2(v) \bigr).
$$

First of all, we need to describe the "(co)adjoint invariants" or, which is the same,  the center of the
corresponding Poisson algebra $S(L)$.
Since $\St(v)$ can be considered as a semisimple Lie algebra in $gl(V)$,  one can use the  polynomial
functions $F_k: L^*(\goth g,\rho, V)=L(\goth g,\rho, V) \to \Bbb
K$   given by
$$ F_k(\Psi) = \Tr_\rho \Psi(v)^k. $$
It is easy to see that  $F_k \in S(L)$, $k=1,2,\dots$

Thus, the commuting polynomials in $S(L)$ constructed by the
argument shift method can be written as follows:
\begin{equation}
F_{\lambda,k} (\Psi) = \Tr_\rho (\Psi(v)+\lambda\Psi_0(v))^k,
\label{invar}
\end{equation}
where $\Psi_0:V\to \goth g$ is a fixed rational section of the stationary
subalgebra fiber bundle (in other words, $\Psi_0 \in L=L(\goth g,\rho, V)$)
satisfying one additional condition: for a generic $v\in V$, the
corresponding element $\Psi_0(v)$ must be regular in $\St(v)$.

The completeness of the set of such polynomials (over $\Bbb K$) is
evident. Indeed, the completeness condition for $L$ is equivalent to the
completeness condition for the functions $\Tr_\rho (X+\lambda A)^k$
defined on $\St(v)$ for generic $v\in V$ (here $X\in\St(v)$ is  variable, $A \in\St(v)$ is
fixed).  But the last condition holds just
because $\St(v)$ is a usual semisimple algebra over $\Bbb C$ (see Theorem \ref{MF}).

%\newpage

\section{
Proof of the Mishchenko-Fomenko conjecture
}

Now we are ready to prove the Mishchenko-Fomenko conjecture.
The following statement reduces the general situation to several separate
cases.

\begin
{lemma} Let $\goth g$ be a Lie algebra over a field $\Bbb K$ of zero
characteristic. Then one of the following statements holds:

(i) $\goth g$ has a commutative ideal $\goth h$ which satisfies at least one of
the two conditions: either  $\dim \goth h >1$  or  $[\goth h, \goth g]\ne 0$;

(ii) $\goth g$ has an ideal $\goth h$ isomorphic to the Heisenberg algebra
$\goth h_m$ and the center of $\goth g$ coincides with the center of
$\goth h$;

(iii) $\goth g=\goth g_0 \oplus \Bbb K$, where $\goth g_0$ is semisimple;

(iv) $\goth g$ is semisimple.
\label{lemma1}
\end{lemma}

 {\it Proof.}
 Consider the radical $\goth r$ of $\goth g$ (if $\goth r$ is trivial, then $\goth g$ is semisimple and we have (iv)).
 Take the chain or ideals:
 $$
\{0\} \subset \goth{r}^{(k)} \subset \goth{r}^{(k-1)} \subset \dots
\subset \goth{r}^{(1)} \subset  \goth{r}^{(0)} = \goth{r}
$$
 where $\goth{r}^{(l+1)} =  [\goth{r}^{(l)}, \goth{r}^{(l)}]$.  Obviously, $\goth{r}^{(k)}$ is a commutative ideal.
 If  $\dim  \goth{r}^{(k)} \ne 1$ or $ \goth{r}^{(k)} $  does not belong to the center $Z(\goth g)$ of $\goth g$, then we get (i).

 Assume that $\dim  \goth{r}^{(k)}=1$ and $\goth{r}^{(k)} \subset Z(\goth g)$.   If the center itself is of dimension greater than 1, then we may take $Z(\goth g)$  as a commutative ideal satisfying (i).

 If $\dim Z(\goth g)=1$, then $\goth{r}^{(k)}$ coincides with   $Z(\goth g)$ and there are  two possibilities:

1) $\goth{r}^{(k)}=\goth{r}$ and then we have case (iii);

2) $\goth{r}^{(k)}$ is contained in the radical $\goth{r}$ as a proper subspace.

In the latter case, consider the ideal $\goth{r}^{(k-1)}$.  If its own center $Z(\goth{r}^{(k-1)})$ is bigger than  $\goth{r}^{(k)}$, then $Z(\goth{r}^{(k-1)})$ is a commutative ideal of dimension greater than 1  and we have case (i).
If  $Z(\goth{r}^{(k-1)})=\goth r^{(k)}$, then $\goth {r}^{(k-1)}$ is a two-step nilpotent Lie algebra
with one-dimensional center, i.e., is isomorphic to the Heisenberg algebra and we have case (iii). $\Box$

It turns out that an induction step (i.e., reducing of dimension) can
naturally  be done in the two first cases (i) and (ii) (see below). In the third and
forth cases no inductive step is needed because a complete commutative subalgebra in $S(\goth g)$
can be constructed by the argument shift method.

%\newpage

Consider the first case (i). Let $\goth h\subset \goth g$ be a commutative ideal.
First of all we give a "differential" description of the polynomials
$f\in\Ann(\goth h)$. For each $x\in \goth g^*$, denote by $h=\pi_{\goth h^*}(x)\in \goth h^*$
its image under the natural projection
$\pi_{\goth h^*}:  \goth g^* \to \goth h^*$. Consider the representation
$(\ad |_\goth h)^*: \goth g \to \mbox{End} (\goth h^*)$ dual to the adjoint one
$\ad|_\goth h: \goth g \to \mbox{End}(\goth h)$ and the corresponding stationary subalgebra
$\St(h)\subset \goth g$ of $h=\pi_{\goth h^*}(x)\in \goth
h^*$.

It is easy to verify the following

\begin{lemma} If $\goth h \subset \goth g$ is an ideal, then $f\in\Ann(\goth h)$ if and only if
$df(x)\in\St(h)$ for any $x\in \goth g^*$.
\label{lemma2}
\end{lemma}

{\it Proof.}
The condition $f\in\Ann(\goth h)$  means that
\begin{equation}
\{ f, \eta\}(x)=\langle x, [df(x),\eta]\rangle = 0 \quad \mbox{for any} \ \ \eta\in\goth h.
\label{cond3}
\end{equation}

Since $\goth h$ is an ideal, this can be rewritten as
$$
0= \langle x, [df(x),\eta]\rangle = \langle h, [df(x),\eta]\rangle =
- \langle  (\ad |_\goth h)^*_{df(x)} h, \eta \rangle,
$$
that is, $(\ad |_\goth h)^*_{df(x)} h=0$, i.e. $df(x)\in\St(h)$, as required. $\Box$

Notice that  (\ref{cond3}) can be rewitten as  $\langle\ad^*_{\goth h} x, df(x)\rangle=0$. In particular, we have
\begin{corollary}
$\St(h)= (\ad^*_{\goth h} x)^\bot = \{\xi\in\goth g~|~\langle\ad^*_{\goth h} x,\xi\rangle=0\}.$
\label{cor1}
\end{corollary}

Since the analysis of differentials is not always an easy task, we give  another version of the above statement, which can be convenient for applications.

\begin{corollary}
Let  $f: \goth g^* \to \Bbb K$  satisfy the condition $f(x +l)=f(x)$ for any $l\in \St(h)$,
$h= \pi_{\goth h^*}(x)$, then
$f \in \Ann(\goth h)$.
\label{cor2}
\end{corollary}

We now describe some "basic" elements in $\Ann(\goth h)$.
Let $\Psi:  \goth h^*
\to \goth g$ be a polynomial map such that $\Psi (h) \in\St(h)$ for
any $h\in \goth h^*$ (among such maps there are, in particular, constant maps
into $\goth h$).  In other words, $\Psi$ is a polynomial section of
the {\it stationary subalgebra fiber bundle} over $\goth h^*$.

The family of such sections is endowed with the natural structure of a Lie
algebra by:
$$
[\Psi_1,\Psi_2](h)=[\Psi_1(h),\Psi_2(h)].
$$

Consider the following polynomial function on $\goth g^*$
\begin{equation}
 f_\Psi (x)= \langle x,\Psi(\pi_{\goth h^*}(x))\rangle.
 \label{func}
\end{equation}

\begin{lemma}
The function $f_\Psi(x)$ belongs to $\Ann(\goth h)$. Moreover, the mapping $\Psi \to f_\Psi$ is a homomorphism of Lie algebras.
\label{lemma3}
\end{lemma}

{\it Proof.} We have
\begin{equation}
df_{\Psi}(x)= d \langle x, \Psi(\pi_{\goth h^*}(x)) \rangle=
\Psi(\pi_{\goth h^*}(x))  +  \langle x, d\Psi(\pi_{\goth h^*}(x)) \rangle.
\label{dfunc}
\end{equation}
The first term $\Psi(\pi_{\goth h^*}(x))=\Psi(h)$ belongs to $\St(h)$ by definition. The second term belongs to $\goth h$, since the section $\Psi$ depends only on
the projection $h=\pi_{\goth h^*}(x)\in \goth h^*$.  Since $\goth h$ is commutative,  we have $\goth h\subset \St(h)$
and, consequently, $df_{\Psi}(x) \in \St(h)$. Thus, $f_\Psi \in
\Ann(\goth h)$ by Lemma \ref{lemma2}.

Furthemore, consider two sections $\Psi_1$ and $\Psi_2$. Denoting $\langle x, d\Psi_i(\pi_{\goth h^*}(x)) \rangle$ by
$\eta_i$, we have
$$
\begin{array}{l}
\{ f_{\Psi_1}, f_{\Psi_2} \}(x)=\langle x,
[\Psi_1(h) + \eta_1, \Psi_2(h)+\eta_2] \rangle = \\
\langle x, [\Psi_1, \Psi_2](h)\rangle +
\langle x, [\Psi_1(h) , \eta_2] \rangle +
\langle x, [\eta_1, \Psi_2(h)] \rangle = \\
f_{[\Psi_1, \Psi_2]} (x) +
\langle h, [\Psi_1(h) , \eta_2] \rangle +
\langle h, [\eta_1, \Psi_2(h)] \rangle
\end{array}
$$

The last two terms vanish since $\Psi_i(h) \in
\St(h)$ and we obtain finally
$$
\{ f_{\Psi_1}, f_{\Psi_2} \}(x)=
f_{[\Psi_1,\Psi_2]}(x).
$$

In other words, the mapping  $\Psi \to f_\Psi$ is a homomorphism of
the algebra of sections into $\Ann(\goth h) \subset
S(\goth g)  $,  as needed. $\Box$

\begin{lemma}
$\trdeg \Ann(\goth h) =\dim \St(h)= \codim \ad^*_{\goth h} x$ for generic $x\in \goth g^*$.
\label{lemma4}
\end{lemma}

{\it Proof.}  The inequality
$$
\trdeg \Ann(\goth h) \le \codim \ad^*_{\goth h} x
$$
is general  and simply means that
"the number of independent invariants cannot be greater than the codimension of a generic orbit".  On the other hand, Lemma \ref{lemma3} explains how one can construct at least  $\dim \St(h)$ algebraically independent  polynomials from $\Ann(\goth h)$,  hence
$$
\trdeg \Ann(\goth h) \ge \dim \St(h).
$$
Finally,  the equality $\dim \St(h)= \codim \ad^*_{\goth h}x$  follows directly from
Corollary \ref{cor1}. $\Box$

This  statement says that $\Ann (\goth h)$ has sufficiently many independent polynomials and
we may apply Proposition \ref{prop2}  (see Remark \ref{rem2}). In other words,   a complete commutative subalgebra in $S(\goth g)$ can be obtained from any two complete commutative subalgebras  $\Cal A_1 \subset  S(\goth h)$ and $\Cal A_2 \subset \Ann(\goth h)$. Also notice that in our case
$S(\goth h)\subset \Ann(\goth h)$ so that we only need to construct a commutative subalgebra $\Cal A$ which is complete in $\Ann(\goth h)$. In other words, we have

\begin{proposition}
Let $\Cal A$ be a complete commutative subalgebra in $\Ann(\goth h)$, then $\Cal A$ is complete in $S(\goth g)$.
\label{prop3}
\end{proposition}

 Another important remark is that $S(\goth h)$ is contained in the center of $\Ann(\goth h)$ so that we may consider  polynomials from $S(\goth h)$ as "new coefficients".   Now we are going to explain how this idea allows us to reduce the problem to a Lie algebra of lower dimension (but over an extended field!).

Let
$p=p(\eta_1,\dots,\eta_l)\in S(\goth h)  $ be an arbitrary polynomial on $\goth h^*$, where
$\eta_1,\dots,\eta_l$ is a certain basis in $\goth h$. If $\Psi: \goth h^* \to \goth g$ is a
polynomial section of the stationary subalgebra fiber bundle, then so is
$p\Psi$.  Besides $[p_1\Psi_1, p_2\Psi_2]=p_1p_2[\Psi_1.\Psi_2]$. This means
that elements from $S(\goth h)$ can be treated as "new coefficients" for the
algebra of sections.  The same is true for $\Ann(\goth h)$:  it is a  module over the ring
$\Bbb K[\goth h^*]=S(\goth h)$ (not only as a commutative algebra of
polynomials but also as a Lie algebra).   Moreover, the
homomorphism of Lie algebras $\Psi \to f_\Psi$ is $\Bbb K[\goth h^*]$-linear.

This observation allows us to pass to a new field of coefficients, namely
$\Bbb K (\goth h^*)=\Frac S(\goth h)  $. To do this correctly we need to
extend all our objects by admitting division by polynomials from
$\Bbb K[\goth h^*]=S(\goth h)$. Instead of $\Ann(\goth h)$ we consider
$$
\Ann_{\rm frac}(\goth h)=\left\{ \left. \frac{f}{g}~\right|~ f\in\Ann(\goth h), g\in S(\goth h)
\right\}
$$

Analogously, instead of polynomial sections $\Psi:\goth h^*\to \goth g$,  we consider
rational ones. As above (see example in Section~3), we
denote the algebra of rational sections by $L(\goth g,(\ad |_\goth h)^*,
\goth h^*)$, and its image in $\Ann_{\rm frac}(\goth h)$ under the mapping $\Psi
\to f_\Psi$ by $L_{\goth h}$.

The crucial point of the proof is that all these objects
$\Ann_{\rm frac}(\goth h)$,  $L(\goth g,(\ad |_\goth h)^*, \goth h^*)$ and $L_{\goth h}$ can
now be treated as Lie algebras over $\Bbb K (\goth h^*)=\Frac S(\goth
h)$. The same is true for the homomorphism  $\Psi \to f_\Psi$.  Moreover, though the Lie
algebra $L_{\goth h}$ is infinite-dimensional over the initial field $\Bbb K$,  it becomes finite-dimensional over $\Bbb K(\goth h^*)$!

\begin
{lemma}  $\dim_{\Bbb K(\goth h^*)  }  L_{\goth h} = \dim_{\Bbb K}\St(h) - \dim_{\Bbb
K} \goth h + 1$, where $\St(h)$ is a generic stationary subalgebra of the representation
$(\ad |_\goth h)^*: \goth g \to \mbox{\rm End} (\goth h^*)$.
\label{lemma6}
\end{lemma}

 {\it Proof.}
To find the dimension of $L_{\goth h}$,  we describe the kernel of the homomorphism
$\Psi
\to f_\Psi$.  It is not hard to see that
$f_\Psi=0$ if and only if  $\Psi(h) \in \Ker (h)$, where
$\Ker(h) \subset \goth h$ is the kernel of the linear functional $h\in \goth h^*$.
The dimension of the subspace of such sections $\Psi$ over $\Frac S(\goth h)  $ is equal
to $\dim \goth h - 1$. Taking into account that the dimension of the algebra of
sections $L(\goth g,(\ad |_\goth h)^*, \goth h^*)$ over $\Bbb K(\goth h^*)$
is equal to the dimension of a generic fiber, i.~e., $\dim_{\Bbb
K}\St(h)$, we immediately obtain the result. $\Box$

\medskip

Thus, we have constructed  a  finite dimensional
subalgebra $L_{\goth h}\subset
\Ann_{\rm frac}(\goth h)$ over the extended field $\Bbb K(\goth h^*)$.
Notice that its dimension is strictly less than $\dim \goth g$ (it concides
with $\dim \goth g$ in the only case, when $\dim \goth h=1$ and simultaneously
$\dim\St(h)=\dim \goth g$, i.e. $\goth h\subset Z(\goth g)$, but exactly this
situation has been excluded from case (i), see Lemma \ref{lemma1}).

Assume that we are able to solve our initial problem (i.e., to construct a complete commutative subalgebra) for the  finite dimensional Lie algebra $L_{\goth h}$ in the sense of the new field $\Bbb K(\goth h^*)$.  It turns out that this leads us immediately to the solution of the problem for  $\goth g$
over the initial field $\Bbb K$.   To see this,  we just need to give some  comments.

    Let ${\Cal A}$ be a complete commutative subalgebra in $S(L_{\goth h})$ in the sense of  $\Bbb K(\goth h^*)$. Without loss of generality,  we shall assume that together with any two polynomials $f$ and $g$ the algebra ${\Cal A}$ contains their product $fg$ and also contains all the constants, i.e. elements from $\Bbb K(\goth h^*)$.

Notice first of all  that  $S(L_{\goth h})$ can naturally be considered as a  subalgebra in $\Ann_{\rm frac}(\goth h)$, since $L_{\goth h} \subset \Ann_{\rm frac}({\goth h})$.  Therefore any commutative subalgebra    ${\Cal A}\subset S(L_{\goth h})$  can be treated as a commutative subalgebra  in $\Ann_{\rm frac}({\goth h})$.

Thus, we can look at ${\Cal A}$ from two different points of view:  either as a subalgebra in $S(L_{\goth h})$ in the sense of the extended field $\Bbb K(\goth h^*)$,  or a subalgebra in $S(L_{\goth h})$ in the sence of the initial field $\Bbb K$ (and then both  $\Cal A$ and $S(L_{\goth h})$ are considered as subalgebras in $\Ann_{\rm frac}({\goth h})$).

We have assumed that  ${\Cal A}$  is complete in $S(L_{\goth h})$ in the  sense of $\Bbb K(\goth h^*)$. Will it be complete in
$S(L_{\goth h})$  in the  sense  of the initial field  $\Bbb K$?  It is not hard to see that the answer is positive.

The next question:  is  this algebra  $\Cal A$ complete in $\Ann_{\rm frac}{\goth h}$)?  The answer is
obviously positive because at a generic point $x\in \goth g^*$, the subspaces  in $\goth g$ generated by the differentials of functions from $S(L_{\goth h})$  and from ${\Cal A}$  are exactly the same (both of them coincide with $\St(h)$, see Lemma \ref{lemma2}).

The last difficulty  is that the functions from ${\Cal A}$ are not polynomial, but rational.
More precisely, they are all of the form  $\dfrac{f}{g}$, where $g \in \Bbb K(\goth h^*)$.  But together with
$\dfrac{f}{g}$ this subalgebra contains both $f$ and $g$ separately.
Therefore, the difficulty can be avoided just by taking the "polynomial" part of ${\Cal A}$,
or simply by multiplying each fraction by its denominator.
After this operation we obtain a certain subalgebra $\Cal A_{\rm pol}$ in $\Ann (\goth h)$ which is obviously commutative and complete (just because the number of independent functions remains the same).
In other words,  after "polynomialization" $\Cal A \to \Cal A_{\rm pol}$ any complete commutative subalgebra ${\Cal A} \subset S(L_{\goth h})$  remains complete in $\Ann(\goth h)$.  Taking into account Proposition \ref{prop3},  we come to the following conclusion.

\begin{proposition}
If  the Mischenko-Fomenko conjecture holds for $L_{\goth h}$   over $\Bbb K(\goth h^*)$, then it holds for  $\goth g$ over the initial field $\Bbb K$.
\label{prop4}
\end{proposition}

Thus,  in case (i)  from Lemma \ref{lemma1},  the problem is reduced to a Lie algebra of smaller dimension.

%%%%%%%%%%%%%%%%%%%%%%%%%%%%%%%%%%%%%%%%%%%%%%%%

Let us now consider the second case. Suppose that algebra $\goth g$ has an ideal
isomorphic to the Heisenberg algebra $\goth h_m$, and the center of $\goth h_m$
coincides with the center of $\goth g$. Recall the structure of the Heisenberg
algebra :
$\goth h_m$ splits into the direct sum of a subspace $V$ of dimension $2m$
and the one-dimensional center $Z(\goth h_m)$ generated by a vector $e$.
For two arbitrary elements  $\xi_1,\xi_2\in V$, their commutator is
defined by $$ [\xi_1,\xi_2]=\omega(\xi_1,\xi_2) e, $$ where $\omega$ is a
symplectic form on $V$.

First we notice several useful properties of
$\goth g$.

 \begin{lemma}
 Let $\goth h_m\subset \goth g$ be an ideal. Then there exists a
subalgebra $\goth{b}\subset \goth g$ such that $\goth g = \goth{b}\oplus  V$
and $\goth{b}\cap \goth h_m=Z(\goth h_m)$. Besides, the subspace
$V\subset \goth{h}_m$ is invariant under the adjoint action of
$\goth b$ and $\goth{b}$ acts on $\goth h_m$ by symplectic
transformations.
\label{lemma7}
\end{lemma}

\noindent {\it Proof.}
We define $\goth{b}$ in the following way:
$$
\goth b = \{ \xi\in\goth g~|~ \ad_\xi (V) \subset V.
$$

Obviously, $\goth b$ is a subalgebra in $\goth g$.  Let us check that any element $\xi\in\goth g$ can be uniquely presented in the form $\xi=\xi_1+\xi_2$, where $\xi_1 \in \goth b$, $\xi_2 \in V$.

 For  $v \in V$, we take
$[\xi,v]\in \goth h_m$ and decompose it with respect to the
subspaces $V$ and $Z(\goth h_m)$:
$$
[\xi,v]=\eta_1+\eta_2, \quad \eta_1 \in V,  \eta_2\in Z(\goth h_m).
$$

Since the center $Z(\goth h_m)$ is one-dimensional $\eta_2$ can be
presented as $\eta_2=l_\xi(v)e$, where $l_\xi:V\to \Bbb K$ is a certain
linear functional. Since $V$ is endowed with a non-degenerate symplectic
structure, this functional can be taken in the form
$l_\xi(v)=\omega(\xi_2,v)$, where $\xi_2 \in V$ is a certain element  which
is uniquely defined by $\xi$. It is easy to see that $\xi-\xi_2$ leaves
the space $V$ invariant:  $$
[\xi-\xi_2,v]=\eta_1+\eta_2-[\xi_2,v]=\eta_1+\omega(\xi_2,v)e-
\omega(\xi_2,v)e=\eta_1\in V.  $$

Thus,  $\goth g=\goth{b}\oplus V$ is a direct sum of the subspaces.
Also it is easy to see that, $\goth{b} \cap \goth h_m=Z(\goth
h_m)$.

We need finally to prove that the representation $\ad:  \goth{b} \to \mbox{End}(V)$
is symplectic, i.e., each transformation
$\ad_\beta:V\to V$ is an element of the symplectic Lie algebra
$sp(V,\omega)$ for any $\beta\in \goth{b}$.

To this end, we use the Jacobi identity. We have:
$$
\ad_\beta
[v_1,v_2]=[\ad_\beta v_1,v_2]+[v_1,\ad_\beta v_2]= \omega(\ad_\beta
v_1,v_2)+\omega(v_1,\ad_\beta v_2).
$$
On the other hand, $[v_1,v_2]$ belongs to the center, therefore $\ad_\beta
[v_1,v_2] =0$.  Thus,
$$
\omega(\ad_\beta v_1,v_2)+\omega(v_1,\ad_\beta
v_2)=0,
$$
which is equivalent to the symplecticity of the represenation
$\ad:  \goth{b} \to gl(V)$. $\Box$

\begin{remark}
{\rm It is not hard to verify that $\ind \goth{b}=\ind \goth g$. The proof is straightforward.
The same result will, however, follow from our consideration below. }
\end{remark}

Following our general idea we need to consider $\goth h_m$ and its annihilator $\Ann(\goth g_m)$.
It turn out that the functions from $\Ann(\goth g_m)$ admit a very natural description.

For any element $\beta\in \goth{b}$ we define a quadratic polynomial
\begin{equation}
f_\beta(x)=\langle \beta, x
\rangle \langle e, x \rangle + \frac{1}{2} \langle \omega^{-1}\bigl(
(\ad_\beta)^*\pi(x) \bigr), x \rangle.
\label{fbeta}
\end{equation}

Here  $\pi : \goth g^* \to V^*$
is the natural projection, $(\ad_\beta)^*:V^*\to V^*$ is the operator dual
to $\ad_\beta: V\to V$,  $\omega$ ---
is a symplectic structure on $V$ treated as a mapping from $V$ to $V^*$
so that $\omega^{-1}$ is an inverse operator from
$V^*$ to $V$, $e$ is a basis element of the center.

\begin{lemma}
$f_\beta \in \Ann(\goth g_m)$.
\label{lemma8}
\end{lemma}

\noindent{\it Proof.} We need to verify the following identity
$$
\langle  x , [df_\beta(x),\eta]\rangle =0
$$
for any  $\eta\in \goth h_m$, $x\in \goth g^*$.

Compute the differential of $f_\beta$. First notice that the
quadratic form  $\langle C x ,
y \rangle=\langle \omega^{-1}\bigl( (\ad_b)^*\pi(x) \bigr), y
\rangle$ is symmetric, therefore  $d\langle
Cx , x \rangle=2Cx$. Hence $$ df_\beta(x) = \beta \langle e,x \rangle +
e \langle \beta, x \rangle + \omega^{-1}\bigl( (\ad_\beta)^*\pi(x)
\bigr).  $$

Then for arbitrary $\eta\in \goth h_m$ we have:
$$
\aligned
&\langle [df_\beta(x),\eta], x \rangle=\\
&\langle [\beta \langle e,x \rangle + e \langle \beta, x \rangle +
\omega^{-1}\bigl( (\ad_\beta)^*\pi(x) \bigr),\eta],x \rangle \\
&\langle e,x \rangle \langle \ad_\beta \eta,x\rangle +
\omega\bigl(\omega^{-1}\bigl( (\ad_\beta)^*\pi(x),\eta) \langle
e,x\rangle = \\ &\langle e,x \rangle \langle \ad_\beta \eta,x\rangle +
\langle (\ad_\beta)^*\pi(x),\eta \rangle) \langle e,x\rangle = \\
&\langle e,x \rangle \langle \ad_\beta \eta,x\rangle - \langle
\pi(x),\ad_\beta \eta\rangle) \langle e,x\rangle = \\ &\langle e,x
\rangle \langle \ad_\beta \eta, x\rangle - \langle x,\ad_\beta
\eta\rangle) \langle e,x\rangle = 0.  \  \Box  \endaligned $$

The next statement is an analog of Lemma \ref{lemma4}.

\begin{lemma}
$\trdeg \Ann({\goth h_m} )= \dim \goth b = \codim \ad^*_{\goth h_m} x$ for generic $x\in \goth g$.
\label{lemma9}
 \end{lemma}

\noindent {\it Proof.}
Here by  $\ad^*$ we denote the coadjoint action of $\goth g$ on $\goth g^*$.  However for the subalgebra $\goth h_m$ we may consider the coadjoint action on its own dual space $\goth h^*_m$.  Denote this action by $\tilde{ad^*}$ for a moment.   Consider two subspaces $\ad^*_{\goth h_m}x$ and $\tilde{\ad}^*_{\goth h_m} h$, where
$x$ is generic in $\goth g^*$ and $h$ is generic in $\goth h_m$.  It is a general and obvious fact that
$$
\dim \ad^*_{\goth h_m}x \ge \dim \tilde{\ad}^*_{\goth h_m} h.
$$
But $\dim \tilde{\ad}^*_{\goth h_m} h = \dim \goth h_m - \ind h_m=2m+1 -1=2m$ so that
$$
\codim \ad^*_{\goth h_m} x  \le \dim \goth g - 2m = \dim \goth b
$$
On the other hand,  Lemma \ref{lemma8} gives us an explicit formula for  $\dim \goth{b}$ independent  polynomials from $\Ann(\goth h_m)$ and, consequently,
$$
 \dim \goth b \le \trdeg \Ann({\goth h_m} )
$$

Taking into account  the general inequality  $\trdeg \Ann(\goth h_m) \le \codim \ad^*_{\goth h_m} x$
we come to the desired conclusion. $\Box$

\medskip

This lemma asserts, in particular, that  $\Ann (\goth h_m)$ has  sufficiently many independent functions so that  we may apply Proposition \ref{prop2}  (see Remark \ref{rem2}). In other words,  we have

\begin{proposition}
Let $\Cal A$ be a complete commutative subalgebra in $\Ann(\goth h_m)$ and $\Cal B$ be a complete commutative subalgebra in $S(\goth h_m)$, then $\Cal A+\Cal B$ is complete in $S(\goth g)$.
\label{prop5}
\end{proposition}

As we see from Lemma \ref{lemma8},    the subalgebra $b$  and the annihilator $\Ann(\goth h_m)$ are closely related.  The  following construction explains this relationship more explicitly.
Instead of $f_\beta$ it will be more convenient to consider the
rational function of the form:
$\tilde f_\beta(x)=f_\beta(x)/\langle e,x\rangle$.

Notice the following remarkable fact which can be verified by a straightforward computation.

\begin{lemma}
The map
$\beta \to \tilde f_\beta$ is  an embedding  (monomorphism) of
$\goth{b}$ into $\text{Frac} (S(\goth g)  )$.
\label{lemma10}
\end{lemma}

The further construction follows the same idea as in the first case (i).  First we need to admit  division by the central elements $g \in S(Z(\goth g))$. Notice  that these elements are just polynomials of one variable $e$,  generator of the center $Z(\goth g)$.  Thus , we consider
$$
\Ann_{\rm frac}(\goth h_m)=\left\{ \left. \frac{f}{g}~\right|~ f\in\Ann(\goth h), g\in S(Z(\goth g)
\right\}
$$

The map $\beta \to \tilde f_\beta$ generates an embedding of $b$ and, consequently, of $S(\goth b)$ into
$\Ann_{\rm frac} (\goth h_m)$.

For applications, it is convenient to rewrite the embedding in dual  terms.
Let $f:\goth b^* \to  \Bbb K$ be a polynomial function on  $\goth b^*$.
Introduce a new function
$\tilde f:\goth g^* \to \Bbb K$ by letting
$$
\tilde f(x) = \tilde f(b+v) = f(b + 1/2  \langle e, b \rangle  \cdot l_v )
$$
where $l_v$ denotes a linear functional on $\goth b$ defined by
$$
l_v(\beta) = {\langle \omega^{-1} ((\ad_\beta)^* v), v\rangle}
$$
and $x=b+v$ is the decomposition dual  to $\goth g = \goth b + V$.

The following statement is just a reformulation of  Lemmas \ref{lemma8} and \ref{lemma10}.

\begin{lemma}
The map  $f \to \tilde f$ is an embedding of $S(\goth b)$ into $\Ann_{\rm frac} (\goth h_m)$.
\end{lemma}

Now it is easy to see that the construction of a complete  commutative subalgebra in $S(\goth g)$ is naturally reduced to the same problem for $S(\goth b)$.

Indeed, suppose we have a complete commutative subalgebra  in $S(\goth b)$.
As before, we assume that this algebra is closed with respect to usual multiplication  and contains
$S(Z(\goth g))$.

 Consider its image
$\tilde{\Cal A}$ in $\Ann_{\rm frac} (\goth h_m)$  under the mapping $f \to \tilde f$.  We claim that $\Cal A$ is complete in $\Ann_{\rm  frac} (\goth h_m)$.  This follows immediately from the fact  that at a generic point,  the subspaces  in $\goth g$ generated by the functions from $\Ann_{\rm frac} (\goth h_m)$  and by the functions of the form $\tilde f$, where $f\in S(\goth b)$ exactly coincide  (since they have the same dimension $\dim \goth b$, see Lemma \ref{lemma9}).  Finally,  to obtain {\it polynomial}  complete commutative subalgebra in $\Ann (\goth h_m)$, we just take the polynomial part $\tilde{\Cal A}_{\rm pol}$ of $\tilde{\Cal A}$,  see above for details.

\begin{proposition}
If  $\goth b$  satisfies the Mischenko-Fomenko conjecture,  then so does $\goth g$.
\label{prop6}
\end{proposition}

Thus, we have shown that in cases (i) and (ii),  the proof of the Mischenko-Fomenko conjecture
can be reduced to the algebra $\goth b$ of  smaller dimension.
The induction argument  completes the proof of Theorem \ref{theo1}.

Notice that the proof is constructive: if we have a complete commutative subalgebra in $S(L_{\goth h})$  or in $S(\goth b)$,  we get  a complete commutative subalgebra in $S(\goth g)$ by using rather simple explicit formulae.

\section{
Examples
}

In this section we show how the above construction works by studying
several examples.  We consider the semidirect sums:

1) $so(n) +_\phi \Bbb R^n$,

2) $sp(2n) + _\phi \Bbb R^{2n}$,

3) $gl(n) + _\phi \Bbb R^n$,
with respect to standard representations.

Recall that  our construction is a step-by-step procedure.  At each step we reduce the dimension
of the Lie algebra under consideration until we come to either one-dimensional  or semisimple Lie algebra.
The first case is the simplest.  After one step we come to a semisimple Lie algebra and then apply the
argument shift method.  The second Lie algebra $sp(2n) +_\phi \Bbb R^{2n}$ needs two steps  (of two different  types corresponding to
cases (1) and (2) from Lemma \ref{lemma1}).
The  affine Lie algebra  $gl(n) +_\phi \Bbb R^n$ is "more complicated":  we never come to the semisimple
algebra, but  have to make $n$ steps  before we finish with the trivial Lie algebra.

We first discuss several general facts.
Consider a semidirect sum $\goth g=\goth k+_\rho  V$ of a Lie algebra $\goth k$ and a commutative ideal $V$.  Its dual space is naturally identified with $\goth k^* + V^*$ and we shall represent elements of $\goth g^*$ as pairs $(M,v)$, where $M\in \goth k^*$, $v\in V^*$.

According to our general approach, we are going to make "reduction"  with respect to  $V$ as a commutative ideal $\goth h$ from Lemma \ref{lemma1}, case 1.
By $\St_{\rho^*} (v)$ we denote the stationary subalgebra of $v\in V^*$ with respect to the dual representation $\rho^*: \goth k \to {\rm End}(V^*)$.
It is easy to see that  the stationary subalgebra  $\St(v)$ considered in Lemma \ref{lemma2} is just the semidirect sum
of $\St_{\rho^*} (v)$ and the ideal $V$. The following
statement  is a reformulation of  Corollary \ref{cor2} in this particular case.

\begin{lemma}
Let $f: \goth g^*  \to \Bbb R$ satisfy the following condition:
\begin{equation}
f(M,v)=f(M+L,v) \quad \mbox{for any} \quad  L\in \St_{\rho^*}(v)^\bot \subset \goth k^*.
\label{cond}
\end{equation}
Then $f\in \Ann (V)$.
\label{lemma11}
\end{lemma}

Condition (\ref{cond}) has a very natural geometrical meaning. Namely, if we think of $v$ as a parameter, then $f(M,v)$  can naturally be considered as a function on $\St_{\rho^*}(v)^*$.
In particular,  this function can be presented in the form
$f(M,v)=f_v (\pi(M))$, where $\pi: \goth k^* \to \St_{\rho^*}(v)^*$ denotes the natural projection.

\begin{lemma}
Let $f(M,u)$ and $g(M,u)$ satisfy {\rm (\ref{cond})}. Then
$$
\{ f(M,v), g(M,v)\} = \{f_v(\pi(M)),  g_v(\pi(M))\}_{\St_{\rho^*} (v)},
$$
where the latter is the Poisson-Lie bracket on $\St_{\rho^*} (v)^*$.
\label{lemma12}
\end{lemma}

The proof of this statement is, if fact, similar to that of Lemma \ref{lemma3} and is based on the simple fact that  $df(M,u)= (X,\eta) \in \goth g$ where  $\eta \in V$, $X\in \St_{\rho^*}(v)\subset \goth k$.

According to the general concept, the construction of a complete commutative subalgebra  in $S(\goth g)$ is reduced to the same problem for $\Ann (V)$. The next statement describes this reduction explicitly.

\begin{lemma}
Consider a set of polynomials  $f_1(M,v), \dots, f_l(M,v)$  satisfying {\rm (\ref{cond})}.  Suppose that  for  generic $v\in V$ they commute as functions on $\St_{\rho^*}(v)^*$ and form a complete commutative
set in $S(\St_{\rho^*}(v))$. Then
$$
\{f_1, \dots, f_l \} \cup V
$$
is a complete commutative set  in $S(\goth g)$.
\label{lemma13}
\end{lemma}

We now pass to the examples.
Consider the Lie algebra $\goth g=e(n)=so(n)+_\phi \Bbb R^n$ (i.e., the Lie algebra of the isometry group of the Euclidean space). The dual space $e(n)^*$ is identified with $e(n)$  by means
of the scalar (non-invariant!) product $\langle (M_1, v_1),(M_2,v_2)
\rangle = \Tr M_1 M_2 +  \langle v_1, v_2\rangle$.

For generic $v\in \Bbb R^n$,  the stationary subalgebra of the standard representation of $so(n)$ is isomorphic to $so(n-1)$.  This stationary subalgebra depends on $v$ as a parameter and is semisimple.
Thus,  a complete commutative set can be constructed by the argument shift method.  According to Lemma \ref{lemma13} we need to construct a set of functions $f_1(M,v), \dots, f_k(M,v)$ such that  for each (generic)
$v$  these functions becomes "the shifts of invariants"  on  the stationary subalgebra of $v$.
As such functions we may consider, for  instance,
$$
f_{\lambda, k}(M,v)=\Tr \bigl(\pr_v (M+\lambda B) \bigr)^k
$$
where  $\pr_v : so(n)=so(n)^* \to \St_{\phi^*} (v)=\St_{\phi^*} (v)^*$ is the orthogonal projection.
It is not hard to see that this projection is given by
$$
\pr_v (M)  = M - \frac{1}{|v|^2}\left( v\otimes (Mv)^\top -
Mv\otimes v^\top\right).
$$

The above functions are not polynomial, but rational.
This problem, however, can easily be avoided by replacing
$\pr_{\St(v)}$ with the map $$|v|^2 \cdot \pr_v: so(n) \to
\St(v)$$
$$
|v|^2 \cdot \pr_v(M)  = |v|^2 M -  v\otimes (Mv)^\bot +
Mv\otimes v^\bot,
$$
which is quadratic in $v$ (and linear in $M$).

As a result we obtain a family of commuting polynomials
$$
\tilde f_{k,\lambda}(M,v) = \Tr \left(|v|^2\pr_v
(M+\lambda B)\right)^k.
$$

The following statement is a particular case of Lemma \ref{lemma13}.  Let $v_i=\langle v, e_i\rangle$ be coordinate linear functions on $\Bbb R^n$ with respect to a certain basis $e_1,\dots, e_n$.

\begin{theorem} {\rm \cite{Ten}}
The functions
$$
v_1, \dots , v_n \ \ \text{and} \ \
\tilde f_{k,\lambda}(M,v), \ \ k=2,4,\dots, [n-1], \ \lambda\in\Bbb R,
$$
generate a complete commutative subalgebra in $S\bigl(e(n)\bigr)$.

\end{theorem}

\begin{remark} {\rm
The above construction was studied by A.S.~Ten in his diploma work \cite{Ten}  two years before Sadetov's proof. In fact, Ten proved this result for any semidirect sum  $\goth k +_\rho V$ if  $\goth k$ is compact. The compactness, however,  can be easily replaced by the assumption that the generic stationary subalgebra of the dual representation $\rho^*: \goth k \to \mbox{End}(V^*)$ is semisimple. }
\end{remark}

The next example is the semidirect product $sp(2n)+_\phi \Bbb R^{2n}$ with respect to the standard
representation.
As above, the elements of $sp(2n) +_\phi \Bbb R^{2n}$  are presented as pairs $(M,u)$, where $M\in sp(2n)$, $u\in \Bbb R^{2n}$.
The dual space is identified with the algebra by
$$
\langle   (M_1,v_1) , (M_2, v_2) \rangle = \Tr M_1M_2 + \Omega(v_1,v_2),
$$
where $\Omega$ is a symplectic form on $\Bbb R^{2n}$.

It is easy to see that the generic stationary subalgebra  $\St_{\phi^*} (v)$  is not semisimple as in the previous case,  but isomorphic to the semidirect sum  $sp(2n-2)+\goth h_{n-1}$, where $\goth h_{n-1}$   is a Heisenberg ideal.
In turn,  $\goth h_{n-1}$ is decomposed into  $(2n-2)$-dimensional symplectic space $V$ and one-dimensional center $Z$.
Such a decomposition is not uniquely defined. To make the choice unique, we choose another element $a\in \Bbb R^{2n}$  such that $\Omega(a,v)\ne 0$.
After this  the subalgebra $sp(2n-2)\subset \St_{\phi^*} (v)$ is defined to be the common stationary subalgebra for $v$ and $a$
\begin{equation}
\St_{\phi^*} (v,a) = \{ A\in sp(2n)~|~  \phi^*(A)a=\phi^*(A)v=0  \},
\label{b}
\end{equation}
the space $V$ is formed by matrices
\begin{equation}
C_\xi=v\otimes(\Omega \xi)^\top + \xi \otimes (\Omega v)^\top
\label{V}
\end{equation}
where $\xi$ belongs to the $(2n-2)$-dimensional subspace
$$\langle v, a\rangle^\Omega=\{ \xi\in \Bbb R^{2n}~|~\Omega(\xi,a)=\Omega(\xi,v)=0\},$$
and the center $Z$ is generated by the matrix
\begin{equation}
C_0=v\otimes (\Omega v)^\top
\label{z}
\end{equation}

Here $\otimes$ denotes  usual matrix multiplication, if we think of $v$ as a column and of $(\Omega v)^\top$ as a row,  at the same time  $\otimes$ is the tensor product of a vector and a covector.

We now apply the general approach to $\St_{\phi^*} (v)=sp(2n-2)+\goth h_{n-1}$ thinking of $v$ as a parameter.
A complete commutative family for $\St_{\phi^*} (v)$ consists of two parts. One is a complete commutative
family for the Heisenberg ideal $\goth h_{n-1}$.
The other is formed by the shifts of $\Ad$-invariants
of $sp(2n-2)$ transmitted into $S(\St_{\phi^*}(v))$ by means of Lemma 10.

The functions corresponding to the Heisenberg ideal are (see (\ref{V}), (\ref{z})):
$$
e(M,v)=\Tr MC_0=\Tr M \, v\otimes (\Omega v)^\top = \Omega(v,Mv)
$$
and
$$
\Tr MC_\xi=\Tr M (v\otimes(\Omega \xi)^\top + \xi \otimes (\Omega v)^\top)=
2\Omega (Mv, \xi)
$$
If we want them to commute, then  $\xi$ must belong to a certain $(n-1)$-dimensional Lagrangian subspace in $\langle v, a \rangle^\Omega$.
For instance, we may take $\xi$ of the form $\xi= \zeta \Omega(u,a)- a \Omega(\zeta,v)$, where $\zeta$ belongs to a certain fixed Lagrangian subspace in $\Bbb R^{2n}$ that contains $a$. In other words, as commuting functions  we can take
$$
f_\zeta(M,v)=\Omega (Mv, \zeta \Omega(v,a)- a \Omega(\zeta,v))=
\left|
\begin{array}{cc}
\Omega( v,\zeta) & \Omega (v,a)\\
\Omega (Mv,\zeta) & \Omega (Mv,a)
\end{array}
\right|
$$

Finally, the shifts of $\Ad$-invariants of $sp(2n-2)=\St_{\phi^*} (v,a)$ take the following form  (after being transmitted into $S(\St_{\phi^*}(v))$ by Lemma \ref{lemma10} and lifted into $S(sp(2n)+_\phi \Bbb R_{2n})$:
$$
f_{k,\lambda}(M,v)=
\Tr \left( \pr_{v,a} \bigl( \Omega(Mv,v) M + Mv\otimes (\Omega Mv)^\top+\lambda B\bigr)\right)^k
$$

 It can be checked that the projection $\pr_{v,a}$ is given by
 $$
 \begin{array}{cl}
 \pr_{v,a}(M)= M - & \Omega(v,a)^{-1}(Ma\otimes (\Omega v)^\top - v \otimes (\Omega Ma)^\top) + \\
& \Omega(v,a)^{-2} \Omega(Ma,a) v \otimes (\Omega v)^\top - \\
& \Omega(v,a)^{-1}(Mv\otimes (\Omega a)^\top - a \otimes (\Omega Mv)^\top) + \\
&\Omega(v,a)^{-2} \Omega(Mv,v) a \otimes (\Omega a)^\top + \\
&\Omega(v,a)^{-2} \Omega(Mv,a) ( a \otimes (\Omega v)^\top +v \otimes (\Omega a)^\top)
\end{array}
 $$

 Thus, to avoid rational functions we  replace $f_{k,\lambda}(M,v)$ by
 $$
\tilde f_{k,\lambda}(M,v)=
\Tr \left(\Omega(v,a)^{2}\pr_{v,a} \bigl( \Omega(Mv,v) M + Mv\otimes (\Omega Mv)^\top+\lambda B\bigr)\right)^k
$$

Here is the final statement.

\begin{theorem} The following functions generate a complete commutative subalgebra
in $S(sp(2n)+_\rho \Bbb R^{2n})$:

1) $v_1,v_2,\dots, v_{2n}$  {\rm(}coordinate functions on $\Bbb R^{2n}${\rm)};

2)$f_\zeta(M,v)$,  where $\zeta$ belongs to a certain
Lagrangian subspace in $\Bbb R^{2n}$ that contains $a$;

3) $e(M,v)$, the function corresponding to the center of  $\St_{\rho^*}(v)$;

4) $\tilde f_{k,\lambda}(M,v)$,  $k=2,4,\dots , 2n$, $\lambda\in \Bbb R$.
\end{theorem}

The last example is the affine Lie algebra $\goth{aff}_n=gl(n,\Bbb R)+\Bbb R^n$.

 Once again we consider $V=\Bbb R^n$ as a commutative ideal and follow our general approach.  The stationary subalgebra of any non-zero element $v\in V^*$ with
respect to the $\Ad^*$-action of $\goth{aff}_n$ on $V^*$ is isomprphic to
$\goth{aff}_{n-1}+\Bbb R^n$, where $\goth{aff}_{n-1} =
\goth{aff}_{n-1}(v) = gl(n-1)+\Bbb R^{n-1}\subset \goth k = gl(n)$ is
the stationary subalgebra of $v$ with respect to the natural action of
$gl(n)$ on $V^*$.
Thus, on the second step of the procedure, we have to deal again with the affine algebra (of smaller dimension) which depends on $u$ as a parameter.  It turns out that repeating this procedure step by step,  we come to the following set of commuting functions.

 Let $\xi_1,\xi_2,\dots,\xi_{n} \subset V=\Bbb
R^n$.  For definiteness, we think of $v$ as a row, and of $\xi_1$ as a
column.  The functions corresponding to the commutative ideal   $\Bbb R^n$ are:
$$ f_{\xi_1} (M,v)= (\xi_1 , v) $$

The functions which correspond to the commutative ideal in the stationary subalgebra $\St(v)=
gl(n-1)+\Bbb R^{n-1}$ take the form
$$
f_{\xi_1,\xi_2}(M,v)= \left| \begin{array}{cc} (\xi_1,v) & (\xi_2,v) \\ (\xi_1,vM) &
(\xi_2, vM) \end{array}\right|
$$

Analogously, on the $k$th step we obtain the functions.
$$
f_{\xi_1,\dots,\xi_k}(M,v)=
\det (a_{ij}),
$$
where  $a_{ij}=(\xi_j,vM^{i-1})$.

\begin{theorem}
The functions $f_{\xi_1,\dots,\xi_k}(M,v)$,  $\xi_i \in \Bbb R^n$,
 $k=1,\dots, n-1$, commute for any values of parameters, i.e.:  $$ \{
f_{\xi_1,\dots,\xi_l},f_{\tilde \xi_1,\dots,\tilde \xi_k}\} =0,
$$ and generate a complete commutative subalgebra in $S(\goth
{aff}_n)$.
\end{theorem}

The proof can be obtained by noticing that if we fix $v$, we obtain the collection of functions on
$\St(v)=\goth{aff}_{n-1}$  just of the same form as the initial functions,  i.e. of the form
$f_{\eta_1,\dots,\eta_k}$ where $\eta_i$ are all orthogonal to the (co)vector  $v$.  It is worth to notice that
$\St(v)=\goth{aff}_{n-1}$ can be naturally interpreted as an affine algebra related to the "orthogonal" complement to $v$, i.e. the subspace  $\{ \eta\in\Bbb R^n~|~ (\eta,v)=0\}$,  $v\in(\Bbb R^n)^*$.
After this remark,  the proof  is obtained by induction.

\section{Two open questions in conclusion}

The Mishchenko-Fomenko conjecture has several natural
generalizations. Actually, the existence of a complete commutative
subalgebra is a very important property to be studied for any
polynomial Poisson algebra. 
One of the most important examples of polynomial Poisson algebra are those of the form
$\Ann(\goth h)$, where  $\goth h$ is a certain subalgebra  of a finite dimensional Lie algebra $\goth g$.

In the particular case of compact Riemannian homogeneous spaces $G/H$,   the existence of  a complete commutative subalgebra in $\Ann(\goth h)$ would  guarantee the integrability of  the geodesic flow on $G/H$ by means of polynomial integrals    (here  $\goth g$ and $\goth h$ are the Lie algebras of  $G$ and $H$ respectively).

\medskip
{\bf Question 1.}  Does $\Ann(\goth h)$ always admit a
complete commutative subalgebra?
\medskip

According to the strong definition of integrability, in addition to
commutativity and completeness of first integrals $f_1, \dots, f_k
\in S(\goth g)$ in the sense Definition 1, one should require the
completeness of each Hamiltonian vector field $X_{f_i}(x)=
\ad^*_{df_i(x)} x$ in the sense that the corresponding  Hamiltonian
flow $\sigma^t_{X_{f_i}}$ is well defined for all $t\in
(-\infty,+\infty)$.

\medskip
{\bf Question 2.}  Consider the complete commutative subalgebra $\mathcal
A \subset \goth g$ constructed in Theorem 1 (see proof in Section
4). Are the Hamiltonian flows of $f\in \mathcal A$ complete?

\end{document}